\documentclass[11pt]{article}

\usepackage[latin1]{inputenc}
\usepackage{epsfig}
\usepackage{latexsym,amssymb,amsmath}
\usepackage{amsthm}
\usepackage{hyperref} 

\newtheorem{thm}{Theorem}[section]

\def\eqref#1{(\ref{#1})}

\newcommand{\todo}[1]
{\vspace{5 mm}\par \noindent \marginpar{\textsc{ToDo}}
\framebox{\begin{minipage}[c]{0.95 \textwidth} \tt #1
\end{minipage}}\vspace{5 mm}\par}

\begin{document}

\title{Explicit equations from orbit reduction: one and two stages}

\author{Viviana Alejandra D\'{\i}az
\\{\footnotesize viviana.diaz@uns.edu.ar}
\\ Departamento de Matem\'atica
\\ Universidad Nacional del Sur,
\\ Av. Alem 1253
\\ 8000 Bahia Blanca, Argentina
\thanks{The authors wish to thank the following institutions for making possible their work on this paper: Universidad
Nacional del Sur (projects PGI and Inter-U);
Universidad Nacional de la Plata (projects PID);
Agencia Nacional de Promoci\'on Cient\'ifica y Tecnol\'ogica (projects PICT);
CONICET (project PIP);
European Community, FP7 (project IRSES).}
\and Marcela Zuccalli
\\{\footnotesize marcezuccalli@gmail.com}
\\ Centro de Matem\'atica de La Plata,
\\ Universidad Nacional de La Plata,
\\ Calle 115 y 47,
\\ 1900 La Plata, Argentina}

\vspace{1cm}
\date{\today}


\maketitle


\begin{abstract}
   It is known that orbit reduction can be performed in one or two stages and it has been proven that the two processes are symplectically equivalent.
   In the context of orbit reduction by one stage we shall write an expression for the reduced two-form in the general case and obtain the
   equations of motion derived from this theory.
    Then we shall develop the same process in the case in which the symmetry group has a normal subgroup to get the reduced symplectic form by two stages
    and the consequent orbit reduced equations.
    In both cases we shall illustrate the method with the example of a rigid body with rotors and compare the obtained equations with the ones given by
    other authors in different frameworks.
\end{abstract}

\tableofcontents

\section{Introduction}

Reduction theory in the presence of symmetry is a very important tool in the study of mechanical systems from the geometric point of
view in relation with the degrees of freedom decrease.
Symmetries are often associated with the concept of integrals of motion. 

There has been a long history of works about the dynamics of mechanical systems using their symmetries. These studies mainly have their origin in the works
of Arnold \cite{arnold66} and Smale \cite{smale70b}. A good exposition of reduction techniques can be found in the book of Abraham and  Marsden
\cite{foundation}.

In the field of classical mechanics, we can deal with reduction theory from two different points of view: Lagrangian and Hamiltonian formalisms, which
are connected by Legendre transform. 

A frame for Lagrangian reduction by stages was developed by Cendra, Marsden and Ratiu in \cite{CMR01a}, presenting Lagrange-Poincar\`e
equations and describing a context in which reduction by stages can be performed. The case in which the configuration space coincides with the symmetry
group
is known as Euler-Poincar{\`e} reduction (see \cite{CHMR98}, \cite{CIM87}, \cite{marsden3}, and \cite{marsdenscheurle}) and, if the group is
abelian, it reduces to the procedure known as Routh reduction (see \cite{marsden3,marsdenscheurle93b,marsdenscheurle}). 
In the framework of Lagrange-Poincar\`e reduction, Cendra and D\'iaz in \cite{CDRedbyStages} generalize the process to the case in which the symmetry group
has a chain of normal subgroups, obtaining an expression for Lie bracket by several stages and the corresponding Lagrange-Poincar\`e equations.

The Hamiltonian version of Lagrange-Poincar\`e reduction was obtained variationally in \cite{CMPR03}. This Hamilton-Poincar\`e theory is connected with
the
one known as Poisson cotangent bundle reduction (see \cite{montgomery86} and \cite{MontMR84}).

For studies in the framework of cotangent bundle reduction we shall mention the referencies \cite{cushmanSniat99}, \cite{HamiltonianReductionAmarillo},
\cite{marsdenPerlmutter00} and \cite{marsdenweinstein74}.

The Dirac procedure of reduction includes Lagrangian, Hamiltonian, and variational points of view simultaneously, and enlarges pre-symplectic and
Poisson reduction. Dirac reduction also encloses the case of
implicit systems and degenerate Lagrangians. In the frame of Dirac theory for symplectic and Poisson cases, we can mention the references
\cite{Courant90,Dorfman93,MR86}.
The case of cotangent bundles of Lie
groups was studied by Marsden and Yoshimura in \cite{YoMa07}, and generalized by the same authors in \cite{YoMa09}. The context for Dirac reduction by
stages was presented in \cite{CMRY09} and \cite{YoMa09}.

A symmetry of a system is usually encoded by a Lie algebra action which is associated with a function from the phase space to the dual of the Lie algebra,
named the momentum map, whose level sets are preserved by the dynamics (see \cite{kostant65,lie1890,marsden3,smale70b,souriau66,souriau69,weinstein83b}).

The structure of the quotient of the preimage of the momentum map of a point 
for compact Lie groups acting over manifolds of finite dimension was developed in  \cite{sjamaarlerman}. If the coadjoint orbit 
is locally closed in the dual of the symmetry group, 
this structure of the quotient was extended for proper actions of Lie groups to the quotient of the preimage of the momentum map of the orbit 
by Bates and Lerman in \cite{bateslerman}. 

The theory of symplectic reduction by stages was originated by the study of semidirect product that was developed in
\cite{guilleminsternberg80,MRW84b,MRW84a,ratiu80tesis,ratiu81,ratiu82}. In \cite{MMPR98} the authors generalize the theory of semidirect product
reduction, considering a symmetry group acting over a symplectic manifold as an extension of one of its normal subgroups, to study reduction by two stages.
A Lagrangian analogous of the reduction for group extensions was done by Cendra,
Marsden and Ratiu in \cite{CMR01a}. 

 In the framework of symplectic reduction, an approach considers the preimage of the momentum map of the coadjoint orbit of a fix point over the symmetry group
 as the reduced space, while the setting known as
 orbit reduction passes to quotient the preimage of the momentum map of the point by the isotropy group.

In this paper, we shall study orbit reduction theory to obtain with purely symplectic techniques an expression in coordinates for the two-form, with the aim of getting in this context associated equations
of motion considering a unique Lie group of symmetry. Then, inspired in the way in which Cendra and
D\'iaz in \cite{CDRedbyStages} have written equations by stages in coordinates in the Lagrangian framework, we study the case in which the symmetry group
has a normal subgroup to get equations in coordinates by two stages in the symplectic orbit approach.

The paper is organized as follows. In section \ref{Orbit reduction}, we shall obtain an expression in coor\-di\-na\-tes for the reduced two-form  and then,
using this expression, we shall write equations in the context of orbit reduction. 
In section \ref{Orbit reduction by two stages} we shall repeat the process in the case in which the symmetry group has a normal subgroup and perform
reduction by two stages, similarly to \cite{TesisMia}. In both cases, we shall illustrate the technique applying it to the case of a rigid body with rotors, getting orbit reduced equations in
one and two stages for this example.

\section{Orbit reduced equations}\label{Orbit reduction}

In this section we shall recall the process of orbit reduction in one stage and write an expression of the symplectic two-form in coordinates. With this
expression, we will obtain the reduced equations of motion in the case in which the system has a symmetry and apply these results to study the
classical example of a rigid body with rotors.

\subsection{Orbit reduction theory}\label{Background}

Now we will gather some results of the theory of orbit reduction related to our work in this paper. The reader can find the proofs and more details in the references \cite{CMR01a,HamiltonianReductionAmarillo,OrtegaRatiuVerde} and \cite{TesisPerlmutter}.

\vspace{.3cm}
Let us consider a connected symplectic manifold $(M, \omega)$ and a free proper canonical (left)
action $\Phi:G\times M\rightarrow M$ of a Lie group $G$ over $(M, \omega).$
We shall note $\Phi(g,m)=\Phi_g(m)=\Phi^m(g),$ indistinctly.
Being $\mathfrak{g}$ the Lie algebra of $G$, each element $\xi \in \mathfrak{g}$ defines a complete vector field $\xi_M$ on $M$ known as the infinitesimal generator associated to $\xi$. 
If $f$ is a differential function on $(M,\omega)$, its associated Hamiltonian vector field $X_{f}$ is defined as the
unique vector field such that $i_{X_{f}} \omega = \mathrm{d}\, f$.

The map $\boldsymbol{J}:M\rightarrow\mathfrak{g}^*$ defined by the relation
 $\langle\boldsymbol{J}(m),\xi\rangle=\boldsymbol{J}^{\xi}(m),$ for all $\xi\in\mathfrak{g}$ and
 $m\in M,$ is called the \textbf{momentum map} of the $G$-action if $\boldsymbol{J}^{\xi}\in C^{\infty}(M)$ is the Hamiltonian function of the infinitesimal generator $\xi_M$. The map $\boldsymbol{J}$ is said \emph{equivariant} if
$\boldsymbol{J}(\Phi(g,m)) = Ad^{*}_{g^{-1}}\boldsymbol{J}(m),$ for all $g \in G$ and $m\in M,$ where $Ad^{*}$ denotes the coadjoint action of
$G$ on $\mathfrak{g}^*$.

\begin{thm}\label{symplectic orbit reduction}
Let $\Phi:G\times M\rightarrow M$ be a free proper canonical
action of the Lie group $G$ on a connected symplectic manifold
$(M,\omega).$ Suppose that this action has an equivariant
associated momentum map \newline
$\boldsymbol{J}:M\rightarrow\mathfrak{g}^*$.  Let
$\mathcal{O}_{\mu}:=G.\mu\subset\mathfrak{g}^*$ be the $G$-orbit
of the point $\mu\in\mathfrak{g}^*$ with respect to the coadjoint
action of $G$ on $\mathfrak{g}^*.$ Then:
\begin{enumerate}
   \item The set $\boldsymbol{J}^{-1}(\mathcal{O}_{\mu})/G$
   is a regular quotient symplectic manifold with the symplectic form
   $\omega_{\mathcal{O}_{\mu}}$ uniquely characterized by the relation
       $\displaystyle i^*_{\mathcal{O}_{\mu}}\omega=\pi^*_{\mathcal{O}_{\mu}}\omega_{\mathcal{O}_{\mu}}+
       \boldsymbol{J}^*_{\mathcal{O}_{\mu}}\omega^+_{\mathcal{O}_{\mu}},$ 
       where $\boldsymbol{J}_{\mathcal{O}_{\mu}}$ is the restriction of $\boldsymbol{J}$ to
       $\boldsymbol{J}^{-1}(\mathcal{O}_{\mu})$ and $\omega^+_{\mathcal{O}_{\mu}}$ is the
       $+$-symplectic structure on the coadjoint orbit $\mathcal{O}_{\mu}$ given by
       $\omega^+_{\mathcal{O}_{\mu}}(\nu)(\xi_{\mathfrak{g}^*}(\nu),\eta_{\mathfrak{g}^*}(\nu))=
       \langle\nu,[\xi,\eta]
       \rangle),$ for arbitrary $\nu\in\mathcal{O}_{\mu}$ and $\xi,\eta\in\mathfrak{g}.$

       The maps $i_{\mathcal{O}_{\mu}}:\boldsymbol{J}^{-1}(\mathcal{O}_{\mu})\hookrightarrow M$ and
       $\pi_{\mathcal{O}_{\mu}}:\boldsymbol{J}^{-1}(\mathcal{O}_{\mu})\rightarrow \boldsymbol{J}^{-1}(\mathcal{O}_{\mu})/G$ are the natural injection and
       projection, respectively. The pair $(\boldsymbol{J}^{-1}(\mathcal{O}_{\mu})/G,\omega_{\mathcal{O}_{\mu}})$ is called the \textbf{symplectic orbit
       reduced space}.
   \item Let $H$ be a $G$-invariant Hamiltonian.
   The Hamiltonian flow $F_t$ of $H$ leaves the connected components of $\boldsymbol{J}^{-1}(\mathcal{O}_{\mu})$ invariant and commutes with the
   $G$-action, so it induces a flow $F_t^{\mathcal{O}_{\mu}}$ on $\boldsymbol{J}^{-1}(\mathcal{O}_{\mu})/G,$ uniquely determined by
       $\displaystyle\pi_{\mathcal{O}_{\mu}}\circ F_t\circ i_{\mathcal{O}_{\mu}}=F_t^{\mathcal{O}_{\mu}}\circ\pi_{\mathcal{O}_{\mu}}.$
    \item The vector field $X_{h_{\mathcal{O}_{\mu}}}$ generated by the flow $F_t^{\mathcal{O}_{\mu}}$ on
    $(\boldsymbol{J}^{-1}(\mathcal{O}_{\mu})/G,\omega_{\mathcal{O}_{\mu}})$ is Hamiltonian with associated
    \textbf{reduced Hamiltonian function} $h_{\mathcal{O}_{\mu}}\in
    C^{\infty}(\boldsymbol{J}^{-1}(\mathcal{O}_{\mu})/G)$ defined by
        $\displaystyle h_{\mathcal{O}_{\mu}}\circ\pi_{\mathcal{O}_{\mu}}=H\circ i_{\mathcal{O}_{\mu}}.$
  \end{enumerate}

\end{thm}

\vspace{.3cm}
Now, let us consider a left smooth action $\Phi:G\times Q\rightarrow Q$ of a Lie group $G$ on a manifold $Q$, and its
associated lifted left action $G\times T^*Q\rightarrow T^*Q$ given
by $g.\alpha_q=T^*_{g.q}\Phi_{g^{-1}}(\alpha_q)$ for $g\in G$ and
$\alpha_q\in T^*_qQ,$ where $T^*_q\Phi_g$ denotes the dual
application of the tangent map of the diffeomorphism
$\Phi_g$. If $\omega$
is the canonical symplectic structure on
$M=T^*Q$, this $G$-action admits an equivariant momentum map
$\boldsymbol{J}: T^*Q\rightarrow\mathfrak{g}^*$ given by
$\langle\boldsymbol{J}(\alpha_q),\xi\rangle=\alpha_q(\xi_Q(q)),$
for all $\xi\in\mathfrak{g}$.

\vspace{.2cm}
It is assumed that the Lie group $G$ acts freely and properly on
the manifold $Q$, and the $G$-action defines a principal bundle
$\pi:Q\rightarrow Q/G.$
Let us consider the adjoint bundle of $Q,$ namely
$\widetilde{\mathfrak{g}}:=(Q\times\mathfrak{g})/G$, which is the
associated vector bundle defined by the adjoint representation of
$G$ over its Lie algebra $\mathfrak{g}$ whose base is $Q/G$ and its fiber is isomorphic
to $\mathfrak{g}.$

From a principal connection $\mathcal{A}$ on the principle bundle
$\pi:Q\rightarrow Q/G$, one can define the vector bundle
isomorphism
$$\alpha_{\mathcal{A}}:TQ/G\rightarrow T(Q/G)\oplus\tilde{\mathfrak{g}}$$
given by $\alpha_{\mathcal{A}}([v_q]):=T_q\pi(v_q)\oplus[q,\mathcal{A}(q)(v_q)],$
where corner brackets denote equivalence classes in the respective quotient manifolds.
Thus, $\displaystyle(\alpha_{\mathcal{A}}^{-1})^*:(T^*Q)/G\rightarrow T^*(Q/G)\oplus\tilde{\mathfrak{g}}^*$
is a diffeomorphism where $\tilde{\mathfrak{g}}^*$ denotes the coadjoint bundle of $Q.$

Given $\mu \in \mathfrak{g}^*$, let
$\mathcal{O_{\mu}}\subset\mathfrak{g}^*$ be the $G$-coadjoint
orbit of $\mu$, it can be verified that
\begin{equation}\label{isomorfprimerespacio}
(\alpha^{-1}_{\mathcal{A}})^*(\boldsymbol{J}^{-1}(\mathcal{O_{\mu}})/G))=T^*(Q/G)
\times_{Q/G}\widetilde{\mathcal{O_{\mu}}}\subset
T^*(Q/G)\oplus\tilde{\mathfrak{g}}^*,
\end{equation}
where
$\widetilde{\mathcal{O_{\mu}}}:=(Q\times\mathcal{O_{\mu}})/G\rightarrow
Q/G$ is the associated fiber bundle to $Q$ considering the coadjoint action of
$G$ on $\mathcal{O_{\mu}}$ and
$$T^*(Q/G)\times_{Q/G}\widetilde{\mathcal{O_{\mu}}}:=\{(\alpha_{[q]},[q,\nu]_G):
\alpha_{[q]}\in T^*_{[q]}(Q/G),q\in Q,\nu\in\mathcal{O_{\mu}}\}.$$

Therefore, the reduced symplectic form
$\omega_{\mathcal{O_{\mu}}}$ of the orbit reduced space
$\boldsymbol{J}^{-1}(\mathcal{O_{\mu}})/G$ (see Theorem \ref{symplectic orbit
reduction}) is
pushed forward by $(\alpha^{-1}_{\mathcal{A}})^*$ to a symplectic
form $\omega_{red}$ on
$T^*(Q/G)\times_{Q/G}\widetilde{\mathcal{O_{\mu}}}.$ This form can be written as (see \cite{HamiltonianReductionAmarillo,TesisPerlmutter})
\begin{equation}\label{forma reducida}
\omega_{red}=\omega_{Q/G}-\beta,
\end{equation}
where $\omega_{Q/G}$
is the canonical symplectic form on the cotangent bundle of the quotient $T^*(Q/G)$ and $\beta$ is the unique two-form on
$\widetilde{\mathcal{O_{\mu}}}$ determined by
$\displaystyle\pi^*_G\beta=\mathrm{d}\boldsymbol{\alpha}+\pi^*_2\omega^+_{\mathcal{O}_{\mu}},$
where
$\pi_G:Q\times\mathcal{O_{\mu}}\rightarrow\widetilde{\mathcal{O_{\mu}}}$
is the projection of the product on the adjoint bundle,
$\pi_2:Q\times\mathcal{O_{\mu}}\rightarrow\mathcal{O_{\mu}}$ is
the projection on the second factor, and $\boldsymbol{\alpha}$ is the one-form on
$Q\times\mathcal{O_{\mu}}$ given by
$\displaystyle\boldsymbol{\alpha}(q,\nu)(v_q,-\text{ad}_{\xi '}^*\nu)=\langle\nu,\mathcal{A}(q)(v_q)\rangle$
for $q\in Q,\ \nu\in\mathcal{O_{\mu}},\ v_q\in T_qQ,\
-\text{ad}^*_{\xi '}\nu\in T_{\nu}\mathcal{O},$ and $\xi
'\in\mathfrak{g}.$ As it has been shown in \cite{TesisPerlmutter}, $\mathrm{d}\boldsymbol{\alpha}$
 has the expression
\begin{equation}\label{diferencialalfa}
\mathrm{d}\boldsymbol{\alpha}(q,\nu)((u_q,-\text{ad}^*_{\xi '}\nu),(v_q,-\text{ad}^*_{\eta '}\nu))=
\langle\nu,[\eta ',\xi]+[\eta,\xi ']+[\xi,\eta]+B(q)(u_q,v_q)\rangle,
\end{equation}
where $B$
is the curvature of $\mathcal{A},$ $q\in Q,\
\nu\in\mathcal{O_{\mu}},\ \xi,\xi ',\eta,\eta ' \in\mathfrak{g},\
u_q,v_q\in T_qQ,$ $\xi$ and $\eta$
are the elements that induce the vertical components of $u_q$ and $v_q$ with respect to the principal connection $\mathcal{A}$.

\subsection{Reduced two-form in coordinates}\label{An expression for the reduced two-form}

\vspace{.3cm}
In this section, we shall obtain a coordinate expression of the 2-form $\omega_{red}$ needed to write equations in
the orbit reduction framework.

\vspace{.3cm}
We consider an element of $T^*(Q/G)\times_{Q/G}\widetilde{\mathcal{O}_{\mu}}$ as a triple $(x,\alpha,\overline{\nu}).$
That is, $x=[q]_G\in Q/G,$ $\alpha\in T^*_x(Q/G),$ $\nu\in \mathcal{O}_{\mu}$ and $\overline{\nu}=[q,\nu]_G\in\widetilde{\mathcal{O}_{\mu}}$ is the adjoint
element of $\nu.$

\vspace{.2cm}
From now on, we shall note $\mathcal{O}$ instead of $\mathcal{O}_{\mu}$ the  coadjoint orbit of a fixed element $\mu$.

Let $(\dot{x}_1,\dot{\alpha}_1,\dot{\overline{\nu}}_1)$ and $(\dot{x}_2,\dot{\alpha}_2,\dot{\overline{\nu}}_2)$ be two vectors belonging to
$T_{(x,\alpha,\overline{\nu})}\left(T^*(Q/G)\times_{Q/G}\widetilde{\mathcal{O}}\right).$ That is,
$(\dot{x}_i,\dot{\alpha}_i)\in T_{(x,\alpha)}T^*(Q/G)$ and
$\displaystyle \dot{\overline{\nu}}_i=\frac{D}{Dt}\overline{\nu}_i\in T_{\overline{\nu}}\widetilde{\mathcal{O}}$
is the covariant derivative of a curve $\overline{\nu}_i$ on $\widetilde{\mathcal{O}},$ for $i=1,2$.

\vspace{.2cm}
To obtain an expression in coordinates for the two-form $\omega_{red},$ since (\ref{forma reducida}) and the fact that $\omega_{Q/G}$ is the canonical form over $T^*(Q/G)$ given by the equation 
\begin{equation*}
   \omega_{Q/G}(x,\alpha)((\dot{x}_1,\dot{\alpha}_1),(\dot{x}_2,\dot{\alpha}_2))=(\mathrm{d}x\wedge\mathrm{d}\alpha)
   ((\dot{x}_1,\dot{\alpha}_1),(\dot{x}_2,\dot{\alpha}_2))=\langle\dot{\alpha}_2,\dot{x}_1\rangle-\langle
   \dot{\alpha}_1,\dot{x}_2\rangle,
\end{equation*}

\noindent we only need to know an explicit expression of $\beta$ which is defined over $\widetilde{\mathcal{O}}$ by
\begin{equation}
\pi^*_G\beta=\mathrm{d}\alpha+\pi_2^*\omega^+_{\mathcal{O}}.
\end{equation}

\vspace{.2cm}
We can see from \cite{CIM87} that if we have a connection $\mathcal{A}$ on a principal bundle
$Q\rightarrow Q/G,$ we are allowed to split a tangent
vector $\dot{\overline{\nu}}\in \widetilde{\mathcal{O}}$
in its horizontal
and vertical part, and write
\begin{equation}
\dot{\overline{\nu}}=[\dot{q}^h,\nu]_G+[q,\dot{\nu}+\mathcal{A}(\dot{q})
\nu]_G,
\end{equation}
where
$\dot{q}^h$ is the horizontal lift of the tangent vector of the curve $[q(t)]_G$ at the point $q\in Q$ and $\dot{\nu}$ is a tangent vector of $\nu\in \mathcal{O}$.

Thus, we can obtain $\beta$ defined over $\widetilde{\mathcal{O}}$ as
$$ \beta\left(\dot{\overline{\nu}}_1,\dot{\overline{\nu}}_2\right)=\beta([\dot{q_1}^h,\nu]_G+[q,\dot{\nu_1}+\mathcal{A}(\dot{q_1})\nu]_G,
[\dot{q_2}^h,\nu]_G+[q,\dot{\nu_2}+\mathcal{A}(\dot{q_2})\nu]_G)=$$
$$\beta(T\pi_G(\dot{q_1}^h,\dot{\nu_1}+\mathcal{A}(\dot{q_1})\nu),T\pi_G(\dot{q_2}^h,\dot{\nu_2}+\mathcal{A}(\dot{q_2})\nu)=
\pi_G^*\beta((\dot{q_1}^h,\dot{\nu_1}+\mathcal{A}(\dot{q_1})\nu),(\dot{q_2}^h,\dot{\nu_2}+\mathcal{A}(\dot{q_2})\nu))=$$
$$(\mathrm{d}\alpha+\pi_2^*\omega^+_{\mathcal{O}})((\dot{q_1}^h,\dot{\nu_1}+\mathcal{A}(\dot{q_1})\nu),(
\dot{q_2}^h,\dot{\nu_2}+\mathcal{A}(\dot{q_2})\nu))=\mathrm{d}\alpha((\dot{q_1}^h,\dot{\nu_1}+\mathcal{A}(\dot{q_1})\nu),
(\dot{q_2}^h,\dot{\nu_2}+\mathcal{A}(\dot{q_2})\nu))+$$
$$\omega^+_{\mathcal{O}}(\dot{\nu_1}+\mathcal{A}(\dot{q_1})\nu,\dot{\nu_2}+\mathcal{A}(\dot{q_2})\nu)$$

\noindent
where $\dot{q}_i^h$ is the horizontal lift of the tangent vector of the curve $[q_i(t)]_G$ at the point $q\in Q$ for $i=1,2.$

\vspace{.2cm}
Since the vertical part of $\dot{q}_1^h$ and $\dot{q}_2^h$ vanish, from (\ref{diferencialalfa}) we have that
$$ \beta\left(\dot{\overline{\nu}}_1,\dot{\overline{\nu}}_2\right)=\langle\nu,B(q)(\dot{q}_1^h,\dot{q}_2^h)\rangle+
\langle\nu,[\varepsilon_1+\mathcal{A}(\dot{q}_1),\varepsilon_2+\mathcal{A}(\dot{q}_2)]\rangle,$$
where $-\text{ad}^*_{\varepsilon_i}\nu=\dot{\nu}_i$ with $i=1,2.$

\vspace{.3cm}
Finally, we can write the reduced two-form as follows:
\begin{equation}\label{formareducidaglobal}
\omega_{red}((\dot{x}_1,\dot{\alpha}_1,\dot{\overline{\nu}}_1),(\dot{x}_2,\dot{\alpha}_2,\dot{\overline{\nu}}_2))=
\langle\dot{\alpha}_2,
\dot{x}_1\rangle-\langle\dot{\alpha}_1,\dot{x}_2\rangle-\langle\nu,B(q)(\dot{q}_{1}^h,\dot{q}_{2}^h)\rangle\\
-\langle\nu,[\varepsilon_1+\mathcal{A}(\dot{q}_1),\varepsilon_2+\mathcal{A}(\dot{q}_2)]\rangle,
\end{equation} where $-\text{ad}^*_{\varepsilon_i}\nu=\dot{\nu}_i$ with $i=1,2.$

\vspace{.3cm}
Note that this expression of the reduced form is the one obtained by Cendra, Marsden, Pekarsky and Ratiu in \cite{CMPR03} but it has been obtained here by symplectic techniques.


\paragraph{Particular cases}
\begin{itemize}
   \item If $Q=G,$ then the fiber bundle $T^*(Q/G)\times_{Q/G}\widetilde{\mathcal{O}}$ can be identified with $\mathcal{O}.$ So it is clear that
\begin{equation*}
\omega_{red}(-\text{ad}^*_{\varepsilon_1}\nu,-\text{ad}^*_{\varepsilon_2}\nu)=-\langle\nu,[\varepsilon_1,\varepsilon_2]\rangle.
\end{equation*}

\vspace{.2cm}
\item If $G$ is abelian, the space $T^*(Q/G)\times_{Q/G}\widetilde{\mathcal{O}}$ can be thought as $T^*(Q/G)$. Thus, the reduced two-form is given by
\begin{equation*}
\omega_{red}((\dot{x}_1,\dot{\alpha}_1),(\dot{x}_2,\dot{\alpha}_2))=\langle\dot{\alpha}_2,\dot{x}_1\rangle-
\langle\dot{\alpha}_1,\dot{x}_2\rangle-\langle\nu,B(q)(\dot{q}_{1}^h,\dot{q}_{2}^h)\rangle.
\end{equation*}
\end{itemize}

\noindent
Both results coincide with the expressions shown in \cite{TesisPerlmutter}.

\subsection{Reduced equations}

Now we shall write the equations of motion in the context of orbit reduction using the expression  (\ref{formareducidaglobal}) obtained in coordinates of the reduced symplectic form.

\vspace{.2cm}
In order to obtain these equations, we shall consider
a $G$-invariant Hamiltonian $H:T^*Q\rightarrow\mathbb{R}$ and its corresponding reduced Hamiltonian $h_{\mathcal{O}}$ (see theorem \ref{symplectic
orbit reduction}) that we will note $h.$ If $X_h\in\mathfrak{X}(T^*(Q/G)\times_{Q/G}\widetilde{\mathcal{O}})$ is the Hamiltonian vector field  associated to the reduced Hamiltonian $h:T^*(Q/G)\times_{Q/G}\widetilde{\mathcal{O}}\rightarrow\mathbb{R},$ then the reduced equations are given by (see \cite{foundation})
\begin{equation*}
\mathrm{i}_{X_h}\omega_{red}=\mathrm{d}h.
\end{equation*}

We can write the Hamiltonian vector field in coordinates as  $X_h=(\dot{x},\dot{\alpha},\dot{\overline{\nu}}).$ So applying the expression (\ref{formareducidaglobal}) for $\omega_{red}$ we
shall obtain the following system of orbit reduced equations

\begin{equation}\label{globalSimplecticRed}
   \left\{\begin{array}{l}
      \displaystyle\frac{\partial h}{\partial x}=-\dot{\alpha}-\langle\nu,B(q)(\dot{q}^h,\ .\
      )\rangle-\langle\nu,[\varepsilon+\mathcal{A}(\dot{q}),.]\rangle  \\
      \\
      \displaystyle\frac{\partial h}{\partial \alpha}=\dot{x} \\
      \\
      \displaystyle\frac{\partial
      h}{\partial\overline{\nu}}=-\langle\nu,[\varepsilon+\mathcal{A}(\dot{q}),.]\rangle=-\text{ad}^*_{\varepsilon+\mathcal{A}(\dot{q})}\nu
   \end{array}
   \right.
\end{equation} where $-\text{ad}^*_{\varepsilon}\nu=\dot{\nu}.$

\paragraph{Orbit reduced equations vs. Hamilton-Poincar\`e equations.}
With variational techniques, Cendra, Marsden, Pekarsky and Ratiu in \cite{CMPR03} have developed a process of reduction for a
Hamiltonian system and obtained the well known Hamilton-Poincar\`e equations which are

\begin{equation}\label{H-P}
   \left\{\begin{array}{l}
      \displaystyle\frac{Dy}{Dt}=-\frac{\partial h}{\partial x}-\langle\overline{\mu},\widetilde{B}(\dot{x},\ .\ )\rangle \\
      \displaystyle\dot{x}=\frac{\partial h}{\partial y} \\
            \\
      \displaystyle\overline{\nu}=\frac{\partial h}{\partial\overline{\mu}}\\
      \displaystyle\frac{D\overline{\mu}}{Dt}=\text{ad}^*_{\overline{\nu}}\overline{\mu}\\
\end{array}
   \right.
\end{equation}
where $\widetilde{B}$ denotes the reduced curvature of the connection.

It can be seen that systems (\ref{globalSimplecticRed}) and (\ref{H-P}) are equivalent. To do that, we may consider the
cannonical identification between elements in the Lie algebra and tangent vectors of the coadjoint orbit, and the fact that $\overline{\mu}$ in (\ref{H-P}) can be identified with the element $\varepsilon+\mathcal{A}(\dot{q})$ in system (\ref{globalSimplecticRed}).

In spite of system (\ref{globalSimplecticRed}) and (\ref{H-P}) are equivalent, it is clear that they are different because they hold in distinct spaces.

Since Lagrange-Poincar\`e equations are equivalent to Hamilton-Poincar\`e system (see \cite{CMPR03}), orbit reduced equations (\ref{globalSimplecticRed}) result
also equivalent to Lagrange-Poincar\`e equations.

\paragraph{Remark.} According to \cite{batessniatycki} and \cite{koonmarsden97b}, the equations (\ref{globalSimplecticRed}) can be easily adapted to the case of nonholonomic systems with horizontal symmetries taking restrictions to the constraint spaces.

\paragraph{Local version.}\label{Orbit reduced equations}

If we consider $Q\rightarrow Q/G$ as a trivial bundle, we have that $\widetilde{\mathcal{O}}$ which is identified with $\mathcal{O}$ is also a trivial bundle, and we can obtain
the following local expression for the symplectic reduced form
\begin{equation*}
\omega_{red}((\dot{x}_1,\dot{\alpha}_1,\dot{\nu}_1),(\dot{x}_2,\dot{\alpha}_2,\dot{\nu}_2))=\langle\dot{\alpha}_2,
\dot{x}_1\rangle-\langle\dot{\alpha}_1,\dot{x}_2\rangle-\langle\nu,B(q)(\dot{q}_{1}^h,\dot{q}_{2}^h)\rangle-
\langle\nu,[\epsilon_1,\epsilon_2]\rangle
\end{equation*}

\noindent where $\epsilon_i\in\mathfrak{g}$ is such that $-\text{ad}_{\epsilon_i}^*\nu=\dot{\nu}_i,$ 
with $i=1,2.$
Having this expression for $\omega_{red},$ we can write
local orbit reduced equations for Hamiltonian systems with symmetry as follows

\begin{equation*}
   \left\{\begin{array}{l}
      \displaystyle\frac{\partial h}{\partial x}=-\dot{\alpha}-\langle\nu,B(q)(\dot{q}^h,\ .\ )\rangle  \\
      \\
      \displaystyle\frac{\partial h}{\partial \alpha}=\dot{x} \\
      \\
      \displaystyle\frac{\partial h}{\partial\nu}=-\langle\nu,[\varepsilon,.]\rangle=-\text{ad}^*_{\varepsilon}\nu 
   \end{array}
   \right.
\end{equation*}

\vspace{.3cm}
\noindent being $\epsilon\in\mathfrak{g}$ the unique element such that $\dot{\nu}=-\text{ad}^*_{\epsilon}\nu.$

\subsection{Rigid body with rotors as an example}

In this section, we are going to study as an illustrative example the case of a rigid body with three rotors in our scheme, with the aim of getting the orbit reduced equations for this system.

\vspace{.2cm}
 Let $Q=SO(3)\times S^1\times S^1\times S^1$ be the configuration space of a rigid body with three rotors and let $G=SO(3)$ be the symmetry group of the system. Let us consider the left action of
 an element $A\in G$ given by $A.(B,\theta_1,\theta_2,\theta_3)=(AB,\theta_1,\theta_2,\theta_3)$ being $(B,\theta_1,\theta_2,\theta_3)$ an arbitrary
 element of $Q.$
It is clear that this action is proper and free. Thus, we have that $Q\rightarrow Q/G=S^1\times S^1\times S^1$ is a trivial bundle which admits a trivial connection.

If $\mu\in\mathfrak{so}(3)^*\equiv\mathbb{R}^3,$ then its $G$-coadjoint orbit $\mathcal{O}$ is a sphere with ratio $\mu.$ Then, the tangent space in a point
$\nu\in\mathcal{O}$ is given by its orthogonal complement. The bundle 
$T^*(Q/G)\times\mathcal{O}$ can be identified with $T^*(S^1\times S^1\times S^1)\times S^2.$ We shall denote $(\dot{\theta}_1,\dot{y}_1,\dot{\nu}_1)$ and
$(\dot{\theta}_2,\dot{y}_2,\dot{\nu}_2)$ two tangent vectors on a given element $(\theta,y,\nu)\in T^*(S^1\times S^1\times S^1)\times S^2.$
So, the reduced two-form $\omega_{red}$ can be written as follows
$$\omega_{red}((\dot{\theta}_1,\dot{y}_1,\dot{\nu}_1),
(\dot{\theta}_2,\dot{y}_2,\dot{\nu}_2))=\langle\dot{y}_2,\dot{\theta}_1\rangle-
\langle\dot{y}_1,\dot{\theta}_2\rangle-\langle\nu,\varepsilon_1\times\varepsilon_2\rangle,$$ where $\epsilon_i=\dot{\nu}_i\times\nu$ for $i=1,2.$

The reduced hamiltonian of this system $h:T^*(S^1\times S^1\times S^1)\times S^2\rightarrow\mathbb{R}$ is given by (see \cite{CMPR03})
$$h(y,\nu)=\frac{1}{2}\sum_{r=1}^3\frac{(\nu_r-y_r)^2}{I_r}+\frac{1}{2}\sum_{r=1}^3\frac{(y_r^2)}{K_r}$$ where $I_r$ and $K_r$ are
the tensors of inertia.

Since the considered connection is trivial, we have that the covariant are usual derivatives and then we can write the orbit reduced equations as follows

\begin{equation*}
\left\{
\begin{array}{rcl}
   \displaystyle \dot{y}_0 &=& 0
   \\
   \dot{\theta}_0^r&=&\displaystyle \frac{y_r-\nu_r}{I_r}+\frac{y_r}{K_r}
   \\
   \displaystyle -\sum_{r=1}^3\nu_r(\epsilon_0\times\epsilon_j)_r &=& \displaystyle\sum_{r=1}^3\frac{\nu_r-y_r}{I_r}
\end{array}
\right.
\end{equation*}
with $1\leq r\leq 3,$ $1\leq j\leq 3,$ and $\nu\in\mathcal{O}=S^2.$

\vspace{.2cm}
These equations coincides with the ones obtained in \cite{CMPR03} with variational techniques but we are focusing in the symplectic setting to get equations by two stages  for this example.

\section{Orbit reduced equations by two stages}\label{Orbit reduction by two stages}

In the case in which the symmetry group of a Hamiltonian system has a normal subgroup we can apply the results of section \ref{An expression for the
reduced two-form} thinking of the normal subgroup as the symmetry group of the system. So we can obtain the symplectic form in which we shall call the first reduced space and the consequent system of orbit reduced equations. After doing that we still have a
remainder symmetry, given by the quotient of the stabilizer (by the big group) of a distinguished coadjoint orbit over the normal subgroup, with which a
second stage of reduction process can be performed to obtain orbit reduced equations by two stages.

\subsection{Setting of orbit reduction by stages}

\vspace{.5cm}
In this subsection, we will recall the principal result of orbit reduction by stages and some isomorphisms between reduced spaces, which will allow us to
take coordinates in order to write the reduced symplectic form by two stages.

\vspace{.2cm}
As in section \ref{Background}, we assume that the Lie group $G$ acts on $Q$ and the action has a cannonical lifting on $(T^*Q, \omega)$ in a free, proper, and
symplectic way.

We also assume that the
group $G$ has a normal subgroup $N$, so $N$ also acts freely, proper and symplectically on $T^*Q$.

Defining $\nu = i^*\mu\in\mathfrak{n}^*,$ where $i:\mathfrak{n}\rightarrow \mathfrak{g}$ is the inclusion of the Lie algebra $\mathfrak{n}$ of the subgroup $N$ in $\mathfrak{g}$, and considering the action of
$N$ on $T^*Q$ instead of the action of $G$ in the theorem \ref{symplectic orbit reduction}, we can think in the momentum map $\textbf{J}_N:T^*Q\rightarrow\mathfrak{n}^*$  and the reduced space $\textbf{J}^{-1}_N(\mathcal{O}_
{\nu})/N \subset T^*Q/N$ where $\mathcal{O}_{\nu}$ denotes the $N$-coadjoint orbit of $\nu.$

Following the techniques shown in \cite{HamiltonianReductionAmarillo}, we shall consider $G_{\mathcal{O}_{\nu}}=\{g\in G:g.\mathcal{O}_{\nu}\subset\mathcal{O}_{\nu}\}$ the
stabilizer of the orbit $\mathcal{O}_{\nu}$ with respect to the $G$-action on $\mathfrak{n}^*$, which is a closed Lie subgroup of $G$
containing $N$ as a normal subgroup.
Then, we have that the Lie subgroup $G_{\mathcal{O}_{\nu}}$ leaves $\textbf{J}^{-1}_N(\mathcal{O}_{\nu})$
invariant and hence, the Lie group $H=G_{\mathcal{O}_{\nu}}/N$ acts on $\textbf{J}^{-1}_N(\mathcal{O}_{\nu})/N$
in a free, proper, and symplectic way. So we can consider the momentum map of the $H$-action on $\textbf{J}_N^{-1}(\mathcal{O}_{\nu})/N$, namely
$\displaystyle\textbf{J}_{\mathcal{O}_{\nu}}:\textbf{J}_N^{-1}(\mathcal{O}_{\nu})/N\rightarrow \mathfrak{h}^*
$ where $\mathfrak{h}^*$ is the dual of the Lie algebra of the group $H$.
For any element $\rho\in\mathfrak{h}^*$ we shall note $\mathcal{O}_{\rho}$ the $H$-coadjoint orbit of $\rho.$

\vspace{.5cm}
In this context the following theorem holds (see \cite{HamiltonianReductionAmarillo}).

\begin{thm}\label{TeoEtapas}
(Orbit Reduction by Stages Theorem)  If the orbit $\mathcal{O}_{\mu}\subset\mathfrak{g}^*$
 satisfies the stages hypothesis and let $\rho\in\mathfrak{h}^*$ defined by the relation
 $\langle\rho,T_e r_{\mathcal{O}_{\nu}}(\xi)\rangle=\langle\mu,\xi\rangle-\langle\tilde{\nu},\xi\rangle,$
for every $\xi\in\mathfrak{g}_{\mathcal{O}_{\nu}}$ and being $r_{\mathcal{O}_{\nu}}:\mathfrak{g}_{\mathcal{O}_{\nu}}\rightarrow
\mathfrak{g}$ the inclusion, where $\mathfrak{g}_{\mathcal{O}_{\nu}}$ denotes the Lie algebra of $G_{\mathcal{O}_{\nu}}$ and $\tilde{\nu}$ is an arbitrary lineal extension of $\nu$ from $\mathfrak{n}^*$ to $\mathfrak{g}^*,$
 then there is a symplectic diffeomorphism
between $\textbf{J}^{-1}(\mathcal{O}_{\mu})/G$ and
$\displaystyle\textbf{J}^{-1}_{\mathcal{O}_{\nu}}(\mathcal{O}_{\rho})/H.$ 
\end{thm}

\vspace{.2cm}
If we consider the principal bundle $Q\rightarrow Q/G_{\mathcal{O}_{\nu}},$ we have a $G_{\mathcal{O}_{\nu}}$-invariant metric on Q and then, we can define
a connection $\mathcal{A}_{G_{\mathcal{O}_{\nu}}}$ given by this metric, see \cite{CMR01a}.
Let $\mathcal{A}_N$ be the $G_{\mathcal{O}_{\nu}}$-invariant connection on the principal bundle $\pi_N:Q\rightarrow Q/N$ associated to the
$G_{\mathcal{O}_{\nu}}$-invariant metric on the principal bundle $Q\rightarrow Q/G_{\mathcal{O}_{\nu}},$ as it is also explained in \cite{CMR01a}. As before (see
(\ref{isomorfprimerespacio})),
we have an isomorphism
$\displaystyle(\alpha_{\mathcal{A}_N}^{-1})^*:\textbf{J}^{-1}_N(\mathcal{O}_ {\nu})/N\rightarrow T^*(Q/N)\times_{Q/N}\widetilde{\mathcal{O}_{\nu}}$ where $\widetilde{\mathcal{O}_{\nu}}=(Q\times\mathcal{O}_{\nu})/N$ is the
adjoint bundle of the $N$-coadjoint orbit of $\nu$. We shall call $T^*(Q/N)\times_{Q/N}\widetilde{\mathcal{O}_{\nu}}$ the \emph{first reduced space}.

Considering the action of $H$ on $Q/N,$ we can choose a connection $A_H$ on the principal bundle $Q/N\rightarrow (Q/N)/H$ such
that, for all $q\in Q$ and all $v\in TQ,$ we have $\mathcal{A}_{G_{\mathcal{O}_{\nu}}}(v)=0$ if and only if $\mathcal{A}_N(v)=0$ and
$\mathcal{A}_H(T\pi_N(v))=0$ (see \cite{CMR01a}). 

\vspace{.2cm}
Now we shall see that $\textbf{J}^{-1}_{\mathcal{O}_{\nu}}(\mathcal{O}_{\rho})/H$ can also be identified with a bundle in
which we can take local coordinates to write the symplectic form by two stages.

\vspace{.3cm}
Since \quad $\textbf{J}^{-1}_{\mathcal{O}_{\nu}} (\mathcal{O}_{\rho})\subseteq\textbf{J}^{-1}_N (\mathcal{O}_{\nu})/N\simeq
T^*(Q/N)\oplus\widetilde{\mathcal{O}_{\nu}},$ we have that
$\displaystyle\textbf{J}^{-1}_{\mathcal{O}_{\nu}} (\mathcal{O}_{\rho})/H
\subseteq(T^*(Q/N)\oplus\widetilde{\mathcal{O}_{\nu}})/H.$ So we can consider the isomorphism induced by the one defined in \cite{CMR01a} for a
Lagrange-Poincar\`e bundle $W$ and an action of a symmetry group $S$ given by
\[ \alpha _{A_S}^{W} := \alpha _{A_S} \oplus \operatorname{id}_{V/H}: T((Q/N) \oplus V)/S \rightarrow T((Q/N)/S) \oplus
\widetilde{\mathfrak{s}} \oplus V/S
\] where
$\alpha _{A_S}$ is the vector bundle isomorphism $\alpha_{A_S} : T((Q/N)/S) \rightarrow T(Q/N)/S \oplus \tilde{\mathfrak{s}}$ and $\tilde{\mathfrak{s}}$ is the adjoint bundle of the Lie algebra $\mathfrak{s}$ of $S.$
Thus, taking $V=\widetilde{\mathcal{O}_{\nu}},$ $W=T(Q/N)\oplus\widetilde{\mathcal{O}_{\nu}},$ and $S=H,$  we obtain an isomorphism
\begin{equation} \label{isodepatas}
\left(\left({\alpha _{A_H}^{W}}\right)^{-1}\right)^* : \left(T^*(Q/N) \oplus \widetilde{\mathcal{O}_{\nu}}\right)/H \rightarrow T^*\left((Q/N)/H\right) \oplus
\widetilde{\mathcal{O}_{\rho}} \oplus (\widetilde{\mathcal{O}_{\nu}}/H)
\end{equation}
being $\widetilde{\mathcal{O}_{\rho}}=(Q/N\times\mathcal{O}_{\rho})/H$ the adjoint bundle of the $H$-coadjoint orbit of $\rho$. 
We will call \newline $T^*((Q/N)/H) \oplus
\widetilde{\mathcal{O}_{\rho}} \oplus (\widetilde{\mathcal{O}_{\nu}}/H)$ \emph{the second reduced space} and then we think the quotient $\textbf{J}^{-1}_{\mathcal{O}_{\nu}} (\mathcal{O}_{\rho})/H$ as included in the second reduced space, that is,  $$\textbf{J}^{-1}_{\mathcal{O}_{\nu}} (\mathcal{O}_{\rho})/H\subseteq T^*((Q/N)/H) \oplus
\widetilde{\mathcal{O}_{\rho}} \oplus (\widetilde{\mathcal{O}_{\nu}}/H).$$

\subsection{Symplectic form by two stages}\label{an expression for the two-form by two stages}

\vspace{.3cm}
Now we will obtain a coordinate version of the symplectic form by two stages, that is, the corresponding
form on the second reduced space.

\vspace{.3cm}
According to the results of the section \ref{An expression for the reduced two-form}, if we think in the symmetry of $N$ and the principal bundle $\pi_N:Q\rightarrow Q/N,$ we can write the reduced two-form in the first reduced space, named $\omega_1,$ as follows.

\vspace{.2cm}
 We shall consider $\nu\in\mathfrak{n}^*$ and $T^*(Q/N)\times_{Q/N}\widetilde{\mathcal{O}}_{\nu}$ with $\mathcal{O}_{\nu}$ the $N$-coadjoint orbit of $\nu,$
 an element $(x,\alpha,\overline{\nu})\in T^*(Q/N)\times_{Q/N}\widetilde{\mathcal{O}}_{\nu},$ i.e., $x=[q]_N\in Q/N,\quad\alpha\in T_x^*(Q/N)$ and
$\overline{\nu}=[q,\nu]_N\in\widetilde{\mathcal{O}_{\nu}},$ and two tangent vectors $(\dot{x}_1,\dot{\alpha}_1,\dot{\overline{\nu}}_1)$ and
$(\dot{x}_2,\dot{\alpha}_2,\dot{\overline{\nu}}_2)$ on the point $(x,\alpha,\overline{\nu}).$ Then, we can write a coordinate expression for the two form
on the first reduced space as (see (\ref{formareducidaglobal}))
\begin{align*}
\omega_1((\dot{x}_1,\dot{\alpha}_1,\dot{\overline{\nu}}_1),(\dot{x}_2,\dot{\alpha}_2,\dot{\overline{\nu}}_2))=\langle\dot{\alpha}_2,
\dot{x}_1\rangle-\langle\dot{\alpha}_1,\dot{x}_2\rangle-\langle\nu,B_N(q)(\dot{q}_{1}^h,\dot{q}_{2}^h)\rangle
-\langle\nu,[\varepsilon_1+\mathcal{A}_N(\dot{q}_1),\varepsilon_2+\mathcal{A}_N(\dot{q}_2)]\rangle,
\end{align*} where $B_N$ is the curvature associated with the connection $\mathcal{A}_N$ on $\pi_N:Q\rightarrow Q/N,$ $\dot{q}_i^h$ is the horizontal lift
of the tangent vector of the curve $[q_i(t)]_N$ in this principal bundle and $\epsilon_i\in\mathfrak{n}$ is such that $\dot{\nu}_i=-\text{ad}^*_{\epsilon_i}\nu,$ for $ 1\leq i\leq 2.$

\vspace{.2cm}
In order to write an explicit expression for the reduced form by two stages, we shall consider a point
$$\displaystyle([x]_H,\gamma,[x,\tau]_H,[[q,\eta]_N]_H)=(y,\gamma,\overline{\tau},[\overline{\eta}]_H)\in T^*((Q/N)/H)\ \underset{(Q/N)/H}{\times}\
\widetilde{\mathcal{O}_{\rho}}\ \underset{Q/N}{\times}\ \widetilde{\mathcal{O}_{\nu}}/H,$$ with
$x\in Q/N,$ $\gamma\in T^*((Q/N)/H),$
 $\tau\in\mathcal{O}_{\rho}\subset\mathfrak{h}^*$ and
$\eta\in\mathcal{O}_{\nu}\subset\mathfrak{n}^*,$ and two tangent vectors
$(\dot{y}_1,\dot{\gamma}_1,\dot{\overline{\tau}}_1,\dot{[\overline{\eta}_1]}_H)\quad$ and
$(\dot{y}_2,\dot{\gamma}_2,\dot{\overline{\tau}}_2,\dot{[\overline{\eta}_2]}_H)\quad$  belonging to
$T_{(y,\gamma,\overline{\tau},[\overline{\eta}]_H)}\left(T^*((Q/N)/H)\ \underset{(Q/N)/H}{\times}\ \widetilde{\mathcal{O}_{\rho}}\ \underset{Q/N}{\times}
\ \widetilde{\mathcal{O}_{\nu}}/H\right).$

\vspace{.5cm}
We will denote $\omega_2$ the form defined over the second reduced space 
$\displaystyle T^*((Q/N)/H)\times\ \widetilde{\mathcal{O}}_{\rho}\times\ \widetilde{\mathcal{O}}_{\nu}/H$ as the image of the induced form $\omega_1$ on the quotient
space via the isomorphism $(({\alpha _{A_H}^{W}})^{-1})^*$ (see (\ref{isodepatas})).

\vspace{.2cm}
Keeping in mind that $\dot{\overline{\tau}}_i\in T_{\overline{\tau}}\widetilde{\mathcal{O}_{\rho}}$ and
$\dot{\overline{\eta}_i}\in T_{\overline{\nu}}\widetilde{\mathcal{O}_{\nu}}$ for $i=1,2$, we can write

\vspace{.2cm}
$\hspace{-.5cm}\displaystyle\omega_2\left((\dot{y}_1,\dot{\gamma}_1,\dot{\overline{\tau}}_1,[\dot{\overline{\eta}_1}]_H),
(\dot{y}_2,\dot{\gamma}_2,\dot{\overline{\tau}}_2,[\dot{\overline{\eta}_2}]_H)\right)=$
\newline
$\displaystyle
\omega_2\left((\dot{y}_1,\dot{\gamma}_1,[\dot{x}^h_1,\dot{\tau}_1+\mathcal{A}_H(\dot{x}_1)\tau]_H,
[[\dot{q}^h_1,\dot{\eta}_1+\mathcal{A}_N(\dot{q}_1)\eta]_N]_H),
(\dot{y}_2,\dot{\gamma}_2,[\dot{x}^h_2,\dot{\tau}_2+\mathcal{A}_H(\dot{x}_2)\tau]_H,
[[\dot{q}^h_2,\dot{\eta}_2+\mathcal{A}_N(\dot{q}_2)\eta]_N]_H)\right),$

\vspace{.2cm}
\noindent
and identify this image with

$$\omega_1\left(\left(\dot{x}_1,T^*_y\pi_H(\dot{\gamma}_1)+\mathcal{A}_H(\dot{x})^*(\dot{\tau}_1+\mathcal{A}_H(\dot{x}_1)\tau),
[\dot{q}^h_1,\dot{\eta}_1+\mathcal{A}_N(\dot{q}_1)\eta]_N\right),\right.$$
$$\left.\left(\dot{x}_2,T^*_y\pi_H(\dot{\gamma}_2)+\mathcal{A}_H(\dot{x})^*(\dot{\tau}_2+\mathcal{A}_H(\dot{x}_2)\tau),
[\dot{q}^h_2,\dot{\eta}_2+\mathcal{A}_N(\dot{q}_2)\eta]_N\right)\right)=$$
$$=\left\langle \dot{x}_1, T^*_y\pi_H(\dot{\gamma}_2)+\mathcal{A}_H(x)^*(\dot{\tau}_2+\mathcal{A}_H(\dot{x}_2)\tau)\right\rangle-
\left\langle \dot{x}_2, T^*_y\pi_H(\dot{\gamma}_1)+\mathcal{A}_H(x)^*(\dot{\tau}_1+\mathcal{A}_H(\dot{x}_1)\tau)\right\rangle$$
$$-\left\langle\eta,B_N(q)(\dot{q}_1^h,\dot{q}_2^h)\right\rangle-\left\langle\eta,
[\varepsilon_1+\mathcal{A}_N(\dot{q}_1),\varepsilon_2+\mathcal{A}_N(\dot{q}_2)]\right\rangle,$$

\vspace{.2cm}
\noindent
where $\varepsilon_i\in\mathfrak{n}$ are such that $\dot{\eta}_i=-\text{ad}^*_{\varepsilon_i}\eta,$ for $1\leq i\leq 2.$

\vspace{.3cm}
Finally, we obtain a coordinate expression for the symplectic form by two stages as follows
\begin{equation}\label{formadosetapas}
\omega_2\left((\dot{y}_1,\dot{\gamma}_1,\dot{\overline{\tau}}_1,[\dot{\overline{\eta}_1}]_H),
(\dot{y}_2,\dot{\gamma}_2,\dot{\overline{\tau}}_2,[\dot{\overline{\eta}_2}]_H)\right)=\left\langle \dot{x}_1,
T^*_y\pi_H(\dot{\gamma}_2)+\mathcal{A}_H(x)^*(\dot{\tau}_2+\mathcal{A}_H(\dot{x}_2)\tau)\right\rangle
\end{equation}
\begin{equation*}
-\left\langle \dot{x}_2, T^*_y\pi_H(\dot{\gamma}_1)+\mathcal{A}_H(x)^*(\dot{\tau}_1+\mathcal{A}_H(\dot{x}_1)\tau)\right\rangle
-\left\langle\eta,B_N(q)(\dot{q}_1^h,\dot{q}_2^h)\right\rangle-\left\langle\eta,
[\epsilon_1+\mathcal{A}_N(\dot{q}_1),\epsilon_2+\mathcal{A}_N(\dot{q}_2)]\right\rangle
\end{equation*}

\vspace{.2cm}
\noindent
where $T^*_x\pi_H$ is a slight abuse of notation to indicate the tangent map of $T^*_x\pi_H$. Also $\mathcal{A}_H(x)^*$ denotes the tangent map of
$\mathcal{A}_H(x)^*$ which is the dual of the
map $\mathcal{A}_H(x):T_x(Q/N)\rightarrow\mathfrak{h}.$

\subsection{Reduced equations by two stages}

In this subsection, we will write orbit reduced equations in coordinates by stages with the expression of the two form obtained above.

We shall consider a given reduced Hamiltonian
$h_2:T^*((Q/N)/H)\times\widetilde{\mathcal{O}}_{\rho}\times\widetilde{\mathcal{O}}_{\nu}/H \rightarrow\mathbb{R}$ and its associated Hamiltonian
vector field $X_{h_2}=(\dot{y},\dot{\gamma},\dot{\overline{\tau}},[\dot{\overline{\eta}}]_H)\in\mathfrak{X}(T^*((Q/N)/H)\times
\widetilde{\mathcal{O}}_{\rho}\times\widetilde{\mathcal{O}}_{\nu}/H).$ Then the equations of motion in the second reduced space are given by the equation
$\displaystyle
\textbf{\text{i}}_{X_{h_2}}\omega_2=\mathrm{d}h_2$ 
that becomes the system of orbit reduced equations by two stages as follows

\begin{equation*}
\left\{
\begin{array}{rcl}
   \displaystyle\frac{\partial h_2}{\partial y}& = &  -T^*_y\pi_H(\dot{\gamma})-\mathcal{A}_H(x)^*(\dot{\tau}+
   \mathcal{A}_H(\dot{x})\tau)-\left\langle\eta,B_N(q)(\dot{q}^{h},.)\right\rangle\\
   \\
   \displaystyle\frac{\partial h_2}{\partial\gamma}&=& \dot{y} \\
   \\
   \displaystyle\frac{\partial h_2}{\partial\overline{\tau}} &=& \mathcal{A}_H(x)(\dot{x})\\
   \\
   \displaystyle\frac{\partial h_2}{\partial[\overline{\eta}]_H} &=&
   -\langle\eta,[\varepsilon+\mathcal{A}_N(\dot{q}),.]\rangle=-\text{ad}^*_{\varepsilon+\mathcal{A}_N(\dot{q})}\eta
\end{array}
\right.
\end{equation*}

\noindent
where $\varepsilon\in\mathfrak{n}$ such that $\dot{\eta}=-\text{ad}^*_{\varepsilon}\eta.$

\vspace{.5cm}
 Equivalently, taking into account that the notation $\mathcal{A}_H(x)^*$ denotes the derivative of the map $\mathcal{A}_H(x)^*$, we can write the term
 $\mathcal{A}_H(x)^*(\dot{\tau})$ as
 $\displaystyle\left\langle\dot{\tau},B_H(x)(\dot{x},.)+\left[\mathcal{A}_H(x)(\dot{x}),\mathcal{A}_H(x)(.)\right]_H\right\rangle$ and finally obtain the
 system
\begin{equation*}
\left\{
\begin{array}{rcl}
   \displaystyle\frac{\partial h_2}{\partial y}& = & \displaystyle
   -T^*_y\pi_H(\dot{\gamma})-\langle\dot{\tau}+\mathcal{A}_H(\dot{x})\tau,B_H(x)(\dot{x},.)\rangle-\left\langle\eta,B_N(q)(\dot{q}^{h},.)\right\rangle\\
   \\
   \displaystyle\frac{\partial h_2}{\partial\gamma}&=& \dot{y} \\
   \\
   \displaystyle\frac{\partial h_2}{\partial\overline{\tau}} &=& \mathcal{A}_H(x)(\dot{x})\\
   \\
   \displaystyle\frac{\partial h_2}{\partial[\overline{\eta}]_H} &=&
   -\langle\eta,[\varepsilon+\mathcal{A}_N(\dot{q}),.]\rangle=-\text{ad}^*_{\varepsilon+\mathcal{A}_N(\dot{q})}\eta
\end{array}
\right.
\end{equation*}

\noindent
where $\varepsilon\in\mathfrak{n}$ such that $\dot{\eta}=-\text{ad}^*_{\varepsilon}\eta.$

\vspace{.3cm}
Let us note that in the case in which
$G_{\mathcal{O}_{\nu}}=G$ this system could be seen as equivalent to the decomposition of Hamilton-Poincar\`e equations, considering the decomposition of the Lie algebra of the symmetry group in the Lie algebra of its normal subgroup and the Lie algebra of the quotient (see the decompositions of the Lie bracket and the principal connection in \cite{CDRedbyStages}).

\paragraph{Local reduced equations by two stages.}\label{Reduced equations in two stages}

If we again consider $Q\rightarrow Q/G$ as a trivial bundle, we get that the coadjoint bundles of the orbits
$\widetilde{\mathcal{O}_{\nu}}$ and $\widetilde{\mathcal{O}_{\rho}}$ are trivials, so we can take a point
$\displaystyle(y,\gamma,\tau,[\nu]_H)\in T^*((Q/N)/H)\times\mathcal{O}_{\rho}\times\mathcal{O}_{\nu}/H$ and two tangent vectors, for $i=1,2,$
\newline
$\displaystyle(\dot{y}_i,\dot{\gamma}_i,\dot{\tau}_i,[\dot{\nu}_i]_H)\in T_{(\dot{y},\dot{\gamma},
\dot{\tau},[\dot{\nu}]_H)}(T^*((Q/N)/H)\times\mathcal{O}_{\rho}\times\mathcal{O}_{\nu}/H).$
Then we can obtain the two stages-\newline
reduced form as follows

\vspace{.2cm}
$\displaystyle\omega_2((\dot{y}_1,\dot{\gamma}_1,\dot{\tau}_1,[\dot{\eta}_1]_H),(\dot{y}_2,\dot{\gamma}_2,\dot{\tau}_2,
[\dot{\eta}_2]_H))=$

\vspace{.2cm}
$\left\langle\dot{x}_1,T^*_x\pi_H(\dot{\gamma}_2)+\mathcal{A}_H(x)^*(\dot{\tau}_2)\right\rangle-\left\langle\dot{x}_2,
T^*_x\pi_H(\dot{\gamma}_1)
+\mathcal{A}_H(x)^*(\dot{\tau}_1)\right\rangle-
\left\langle\eta,B_N(q)(\dot{q}_1^h,\dot{q}_2^h)\right\rangle-\langle\eta,[\epsilon_1,\epsilon_2]\rangle.$

\vspace{.5cm}
So, in this case we get the following system of local orbit reduced equations by two stages

\begin{equation*}\label{ecuac2etapas}
\left\{
\begin{array}{rcl}
   \displaystyle\frac{\partial h_2}{\partial y}& = &
   -T^*_x\pi_H(\dot{\gamma})-\mathcal{A}_H(x)^*(\dot{\tau})-\left\langle\eta,B_N(q)(\dot{q}^{h},.)\right\rangle\\
   \\
   \displaystyle\frac{\partial h_2}{\partial\gamma}&=& \dot{y} \\
   \\
   \displaystyle\frac{\partial h_2}{\partial\tau} &=& \mathcal{A}_H(x)(\dot{x})\\
   \\
   \displaystyle\frac{\partial h_2}{\partial[\eta]_H} &=& -\langle\eta,[\epsilon,.]\rangle=-\text{ad}^*_{\epsilon}\eta
\end{array}
\right.
\end{equation*}

 Remembering again the notation $\mathcal{A}_H(x)^*$ for the derivative of the map $\mathcal{A}_H(x)^*,$ we can write the system as 

\begin{equation*}
\left\{
\begin{array}{rcl}
   \displaystyle\frac{\partial h_2}{\partial y}& = & -T^*_x\pi_H(\dot{\gamma})-\left\langle\dot{\tau},B_H(x)(\dot{x},.)\right\rangle
   -\left\langle\eta,B_N(q)(\dot{q}^{h},.)\right\rangle\\
   \\
   \displaystyle\frac{\partial h_2}{\partial\gamma}&=& \dot{y} \\
   \\
   \displaystyle\frac{\partial h_2}{\partial\tau} &=& \mathcal{A}_H(x)(\dot{x})\\
   \\
   \displaystyle\frac{\partial h_2}{\partial[\eta]_H} &=& -\langle\eta,[\epsilon,.]\rangle=-\text{ad}^*_{\epsilon}\eta
\end{array}
\right.
\end{equation*}

\vspace{.2cm}
This system of equations is in correspondence via Legendre transform to local Lagrange-Poincar\`e equations by two stages obtained in \cite{CDRedbyStages} in the comparable case, that is, $G_{\mathcal{O}_{\nu}}=G.$

\subsection{Rigid body with rotors revisited}

In this subsection, we shall again consider the rigid body with three rotors and the same configuration space as before but a different group of symmetry, to write equations of motion in the context of orbit reduction by two stages.

Let $Q=SO(3)\times S^1\times S^1\times S^1$ and $G=Q$ acting by the left action given by $(B,\alpha_1,\alpha_2,\alpha_3)(A,\theta_1,\theta_2,\theta_3)=(B.A,\alpha_1+\theta_1,\alpha_2+\theta_2,\alpha_3+\theta_3).$  We shall consider $N=SO(3)$ as a normal subgroup of $G$
and the $N$-coadjoint orbit of an element $\nu\in\mathfrak{n}^*$ denoted $\mathcal{O}_{\nu}\subset\mathfrak{n}^*.$ Thus, we have that
$\mathcal{O}_{\nu}\simeq S^2.$
Then, we can write $\widetilde{\mathcal{O}_{\nu}}=(Q\times\mathcal{O}_{\nu})/N=(SO(3)\times (S^1)^3\times S^2)/SO(3)=(S^1)^3\times S^2.$
In this case, we can obtain that the stabilizer $G_{\mathcal{O}_{\nu}}$ is the group $G$ and so $H=G_{\mathcal{O}_{\nu}}/N=G/N=(S^1)^3.$

If we consider an element $\rho\in\mathcal{O}_{\nu}\simeq S^2$ satisfying the condition in the theorem \ref{TeoEtapas},
since the action of $H$ on $Q/N=(S^1)^3$ leaves it invariant, we have that the $H$-coadjoint orbit of $\rho$ can be identified with the element $\rho$ itself.
Thus, $\widetilde{\mathcal{O}_{\rho}}=
(G/N\times\mathcal{O}_{\rho})/H=((S^1)^3\times\{\rho\})/(S^1)^3\equiv\mathcal{O}_{\rho}=\{\rho\}.$

\vspace{.2cm}
The $H$-invariant reduced Hamiltonian of the system $h_1:T^*(S^1)^3\times\widetilde{S^2}\rightarrow\mathbb{R}$ is given by (see \cite{CMPR03})
$$h_1(y,\nu)=\frac{1}{2}\sum_{r=1}^3\frac{(\nu_r-y_r)^2}{I_r}+\frac{1}{2}\sum_{r=1}^3\frac{y_r^2}{K_r}.$$

We can consider the second stage of reduction of the Hamiltonian, that is, the reduction of $h_1$ by $H$ that we shall note $h_2:\left(T^*(S^1)^3\times\widetilde{S^2}\right)/(S^1)^3\simeq\widetilde{S^2}\rightarrow\mathbb{R}$ given by
$$h_2(\nu)=h_1(0,\nu)=\frac{1}{2}\sum_{r=1}^{3}\frac{\nu_r^2}{I_r}.$$

On the other hand, we know that
$T\left(T^*\left((Q/N)/H\right)\times\widetilde{\mathcal{O}}_{\rho}\times\widetilde{\mathcal{O}}_{\nu}/H\right)=
T\left(T^*\left((S^1)^3/(S^1)^3\right)\times\{\rho\}\times S^2\right)\simeq TS^2.$
So, $\omega_2$ is a 2-form over $S^2.$ Then if $\nu\in S^2,$ from formula (\ref{formadosetapas}) for this particular case, we have \quad $\omega_2(\nu)(\dot{\nu_1},\dot{\nu_2})=-\langle\nu,[\varepsilon_1,\varepsilon_
2]\rangle$ where $\dot{\nu}_i=\text{ad}^*_{\varepsilon_i}\nu,$ for $1\leq i\leq 2.$
In other words,
$\omega_2(\nu)(\dot{\nu_1},\dot{\nu_2})=-\langle\nu,\varepsilon_1 \times\varepsilon_ 2\rangle$ where
$\varepsilon_i=\dot{\nu_i}\times\nu,$ for $1\leq i\leq 2.$
Then, the orbit reduced equation by two stages for this example is given by $$-\text{ad}^*_{\epsilon}\nu=\dot{\nu}=\nu\times\epsilon,$$ or equivalently

\begin{equation*}
\left\{
\begin{array}{rcl}
   \displaystyle \frac{\nu_1}{I_1}&=&pr_1(\nu\times\epsilon)
   \\
   \\
    \displaystyle \frac{\nu_2}{I_2}&=&pr_2(\nu\times\epsilon)
   \\
   \\
    \displaystyle \frac{\nu_3}{I_3}&=&pr_3(\nu\times\epsilon)
\end{array}
\right.
\end{equation*}

\vspace{.2cm}
\noindent
where $pr_i$ is the projection on the $i$-th factor.

\paragraph{Posible future research.}

This work could naturally suggest further investigations in different directions. In the first place we can mention the study of other examples of interest in mechanics to write orbit reduced equations, process that can give additional information about the systems under consideration. Another interesting topics are  the generalization of these techniques to the case of non-equivariant momentum map by considering its associated cocycle, and to the situation in which the symmetry group has a chain of normal subgroups leading to a theory of orbit reduction by several stages.

\paragraph{Acknowledgements.} We would like to thank Eduardo Garc\'ia-Tora\~no Andr\'es and Jorge Solom\'in for the useful discussions.

\bibliographystyle{plain}
\bibliography{bibliogtesis}

\end{document}